\documentclass[%
reprint,
superscriptaddress,
amsmath,amssymb, 
aps,
article,prr,
floatfix,
]{revtex4-2}

\usepackage{graphicx}
\usepackage[T1]{fontenc}
\usepackage{amsfonts}
\usepackage{hyperref}
\usepackage{color}
\usepackage{cancel}
\usepackage{bm}
\usepackage{enumitem}
\usepackage{soul}

\setlist[enumerate]{leftmargin=*}

\newcommand{\LC}{\ell c}
\newcommand{\sgn}{\operatorname{sgn}}

\newcommand{\Kctwo}{K_{c2}}

\begin{document}

\title{Time-optimal synchronisation to self-sustained oscillations under bounded control}
\author{C.~Ríos-Monje}
 \affiliation{Física Teórica, Multidisciplinary Unit for Energy Science, Universidad de Sevilla, Apartado de
  Correos 1065, E-41080 Sevilla, Spain}
\author{C. A. Plata}
 \affiliation{Física Teórica, Multidisciplinary Unit for Energy Science, Universidad de Sevilla, Apartado de
  Correos 1065, E-41080 Sevilla, Spain}
\author{D. Guéry-Odelin}
\affiliation{Laboratoire de Collisions Agr\'egats R\'eactivit\'e, Universit\'e de Toulouse, 
CNRS, 118 Route de
Narbonne, 31062 Toulouse CEDEX 4, France}
\author{A.~Prados} 
 \affiliation{Física Teórica, Multidisciplinary Unit for Energy Science, Universidad de Sevilla, Apartado de
  Correos 1065, E-41080 Sevilla, Spain}
\date{\today}

\begin{abstract}
Incorporating force bounds is crucial for realistic control implementations in physical systems. Here, we investigate the fastest possible synchronisation of a Liénard system to its limit cycle using a bounded external force. To tackle this challenging non-linear optimal control problem, our approach involves applying Pontryagin's Maximum Principle with a combination of analytical and numerical tools. We show that the optimal control develops a remarkably complex structure in phase space as the force bound is lowered. Trajectories rewound from the limit cycle's extreme points turn out to play a key role in determining the maximum number of control bangs for optimal connection. We illustrate these intricate features using the paradigmatic van der Pol oscillator model.
\end{abstract}

\maketitle

\section{Introduction}

Self-sustained oscillations are a type of oscillatory behaviour that persists over time without any external driving~\cite{ginoux_van_2012,jenkins_self-oscillation_2013}. This phenomenon is ubiquitous in nature, playing a fundamental role in many fields: active membranes~\cite{gomez-marin_self-sustained_2007}, climate change~\cite{arnaut_self-sustained_2020}, cancer cell dynamics~\cite{bruckner_stochastic_2019}, circadian rhythms~\cite{leloup_limit_1999,roenneberg_modelling_2008,seki_evolution_2022}, electronic resonators~\cite{hu_phototunable_2021,zhang_frequency_2024,zhang_slow_2025}, network dynamics~\cite{kawai_reservoir_2023,liu_recovery_2025}, or quantum devices~\cite{tabanera-bravo_stability_2024}.  Generally, self-sustained oscillations are only achieved in the infinite time limit, appearing as stable limit cycles in the phase space of the system~\cite{he_determination_2003}. In this context, the Liénard equation describes a broad family of non-linear oscillators~\cite{strogatz_nonlinear_2024}, including van der Pol's~\cite{van_der_pol_lxxxviii_1926}, that exhibit such self-sustained oscillations.

The field of shortcuts, first developed for quantum systems~\cite{guery-odelin_shortcuts_2019} and later transposed to classical systems, both with deterministic and stochastic dynamics~\cite{guery-odelin_driving_2023}, aims to engineer protocols that accelerate the connection between given initial and target states. Very recently, this idea has been employed to investigate synchronisation to a limit cycle in a finite time by applying an external force: first, devising non-optimal connections~\cite{impens_shortcut_2023}; later, protocols minimising the energetic cost~\cite{rios-monje_optimal_2024}. Minimum work connections involve Dirac-delta peaks of the force at the initial and final times~\cite{rios-monje_optimal_2024}, similarly to the situation found in stochastic thermodynamics for underdamped systems~\cite{gomez-marin_optimal_2008,faure_active_2025,baldovin_optimal_2025}, when arbitrarily large values of the driving force are possible.

In any actual physical system, the applied driving $F(t)$ force is bounded---in the simplest scenario, by the non-holonomic constraint $|F(t)|<K$. Aside from limiting the practical applicability of protocols with impulsive forces, like those in Refs.~\cite{gomez-marin_optimal_2008,rios-monje_optimal_2024}, the force bound entails that a minimum time is expected to emerge for the protocol. Driving the system as fast as possible to its long-time state is especially relevant for small dissipation, since the natural relaxation time becomes very long.

In the context of quantum shortcuts, minimum time connections are related to the quantum speed limit and the quantum brachistochrone problem, i.e.~finding the protocol that connects given initial and final states in the shortest time~\cite{deffner_quantum_2017,boscain_introduction_2021,dionis_time-optimal_2023,fanchiotti_quantum_2024,morandi_optimal_2025}. Similar minimum time problems has been recently explored in many fields, e.g.~\cite{sun_exploring_2022,fernandez-cara_analysis_2025}; in the context of classical shortcuts, minimum time protocols have been analysed for connecting equilibrium states and also non-equilibrium steady states~\cite{plata_finite-time_2020,patron_thermal_2022,patron_minimum_2024,baldovin_optimal_2025}. Yet, optimising the connection towards a time-dependent state in the presence of a non-holonomic constraint remains a challenge.

We address the above challenge by optimising, in terms of connection time, the synchronisation of the general class of non-linear oscillators described by the Liénard equation to its limit cycle, by applying a bounded external force. Optimal control theory is the mathematical framework for this problem. Specifically, we are dealing with a minimum time problem, which can be solved using Pontryagin's Maximum Principle (PMP)~\cite{pontryagin_mathematical_1987,liberzon_calculus_2012,boscain_introduction_2021}. PMP provides necessary conditions for the optimal control, i.e.~the optimal driving force that minimises the connection time to the limit cycle and the corresponding optimal trajectory of the system in phase space.

After presenting the Liénard equation, we seek the fastest synchronisation to its limit cycle. PMP allows us to derive key analytical results and find the general structure of the optimal driving: the optimal protocol is always bang-bang, alternating time windows with $F(t)=\pm K$. Notably, a very complex optimal control structure emerges as the bound decreases, with an increasing number of bangs. The analytical results are derived for the general Liénard equation, and are  illustrated in the paradigmatic case of van der Pol's oscillator. 

\section{Model}

The Liénard equation~\cite{chandrasekar_unusual_2005,iacono_class_2011,messias_time-periodic_2011,ghosh_lienard-type_2014,shah_conservative_2015,turner_maximum_2015,gine_lienard_2017,kingston_extreme_2017,suresh_influence_2018,mishra_chimeras_2023,strogatz_nonlinear_2024} describes a wide family of non-linear oscillators with potential $V(x)$ and non-linear damping force $-\mu h(x)\dot{x}$:
\begin{equation}
\label{eq:lienard}
\ddot{x} + \mu h(x) \dot{x} + V'(x) = 0,
\end{equation}
where \(\mu\) is a constant called damping coefficient, and  $h(x)\in\mathcal{C}^1$ and $V(x)\in\mathcal{C}^2$ are even functions. Equation~\eqref{eq:lienard} is written in dimensionless variables. It has a unique and stable limit cycle, denoted by \(\chi_{\LC}(x,\dot{x};\mu) = 0\), under quite general conditions~\footnote{Specifically, the conditions are~\cite{strogatz_nonlinear_2024}: 
(i) \(V(x)\) is a confining potential with only one minimum at \(x=0\), and (ii) the function  $\xi(x) = \int_{0}^{x} h(x') dx'$
has the properties: (a) \(\xi(x)\) has only one positive zero at \(x=a\), with \(\xi(x)<0\) for \(0<x<a\) and \(\xi(x)>0\) for \(x>a\), and (b) \(\xi(x)\) is non-decreasing for \(x>a\)
, with \(\lim_{x\to\infty}\xi(x)=+\infty\).}.

Our goal is to synchronise this system, i.e.~to drive it from an arbitrary initial condition to its self sustained oscillations in a finite time by applying an external driving force \(F(t)\). 
Thus, we add $F(t)$ to the right-hand-side (rhs) of Eq.~\eqref{eq:lienard}. 
The resulting driven Liénard equation is conveniently rewritten in phase space, introducing  $\bm{x}^{T}(t)\equiv (x_1(t),x_2(t))$, $x_1 = x$, $x_2 = \dot{x}$, as
\begin{equation}
     \dot{\bm{x}} \equiv \begin{pmatrix}
      \dot{x}_1 \\ \dot{x}_2
    \end{pmatrix}= \bm{f}(\bm{x},F)\equiv \begin{pmatrix}
    x_2 \\ -\mu h(x_1)x_2-V'(x_1)+F
    \end{pmatrix}.
\label{eq:lienard_driven_vector}
\end{equation}

Among all possible solutions of the synchronisation problem, we are interested in that minimising the time taken to reach the limit cycle when the force is bounded, 
\begin{equation}
    |F(t)|\leq K, \quad \forall t>0.
\label{eq:bound_force}
\end{equation}
Therefore, the optimisation problem to be solved is 
\begin{equation}
  \min_{x\in\mathcal{C}^1,\;|F(t)|\leq K} t_f = \int_{0}^{t_f} dt,
\label{eq:min_prob}
\end{equation}
where the system dynamics verifies Eq.~\eqref{eq:lienard_driven_vector}, the initial phase-space point is fixed, and the final phase-space point \(\bm{x}_f\equiv \bm{x}(t_f)\) belongs to the limit cycle,
\begin{equation}
    \bm{x}^{T}(0) = (x_{10},x_{20}), \quad \chi_{\LC}(\bm{x}_{f};\mu)=0.
\label{eq:initial-and-final-cond}
\end{equation}
The limit cycle splits phase space into two regions: (i) interior points, chosen to be defined by \(\chi_{\LC}(\bm{x};\mu)<0\) and (ii) exterior points, where \(\chi_{\LC}(\bm{x};\mu)>0\). 

PMP provides the necessary conditions for determining the optimal control $F^*(t)$ (here, the driving force) that minimises the connection time, along with the corresponding optimal trajectory $\bm{x}^*(t)$ over the time interval $t \in [0, t_f]$. 
At its core, PMP relies on the definition of the Pontryagin Hamiltonian $\mathcal{H}$, given by:
\begin{align}
   \mathcal{H}(\bm{x},\bm{p},p_0;F)  & =   \bm{p}\cdot\bm{f}(\bm{x},F) + p_0 \nonumber \\  
   & =   p_1 f_1(\bm{x},F) + p_2 f_2(\bm{x},F)+ p_0. 
   \label{eq:hamiltonian}
\end{align}
According to PMP, there exists a costate vector \(\bm{p}^*(t)\) and a constant \(p_0^* \le 0\), with \(\{\bm{p}^*,p_0^*\}\neq \{\bm{0},0\}\) for all \(t\in[0,t_f]\)~\footnote{The constant $p_0^*\le 0$ is sometimes called the abnormal multiplier and it is a key ingredient of PMP. In non-degenerate cases, $p_0\ne 0$ and it is thus typically normalised to $p_0^*=-1$; in some degenerate cases, $p_0^*=0$. See Chapter 4 of Liberzon's book~\cite{liberzon_calculus_2012} for more details.}. These satisfy the following three conditions, 
(i) canonical system equations: the time evolution of $\{\bm{x}(t),\bm{p}(t)\}$ is governed by
\begin{equation}
  \dot{\bm{x}}^* = \nabla_{\bm{p}^*}\mathcal{H}^*, \quad \dot{\bm{p}}^* = - \nabla_{\bm{x}^*}\mathcal{H}^*,
\label{eq:canon_eq}
\end{equation}
together with the boundary conditions given in Eq.~\eqref{eq:initial-and-final-cond}; (ii) Hamiltonian maximisation condition: $\forall t\in[0,t_f]$ and  $\forall F\in[-K,+K]$, $\mathcal{H}$ attains its maximum value, which is zero,  as a function of the driving force at \(F^*\),
\begin{equation}
    \mathcal{H}(\bm{x}^*,\bm{p}^*,p_0^*;F)\leq 0=\mathcal{H}^*\equiv \mathcal{H}(\bm{x}^*,\bm{p}^*,p_0^*;F^*),
    \label{eq:hamiltonian_optimal}
\end{equation}
and
(iii) transversality condition: at the final time \(t_f\), the costate vector is orthogonal to the tangent space of the limit cycle,
\begin{equation}
  \bm{p}_f^* \cdot \bm{d}= 0,\quad \forall\bm{d}\in T_{\bm{x}^*_f}\{\chi_{\LC}=0\}.
\label{eq:transversality}
\end{equation}
Asterisks are dropped hereafter, as all variables refer to the optimal solution.

PMP provides a framework for determining the optimal driving force and solving the minimum time synchronisation problem. However, the non-linearity of the involved equations prevent us from obtaining a full analytical solution for the optimal protocol. Yet, our analysis of PMP makes it possible to discern the structure of the optimal driving force, also enabling us to obtain a complete solution of the optimisation problem numerically.

\section{Main results}

We now present the main results for the minimum time problem stemming from PMP; detailed derivations are in the Appendix. Remarkably, all results apply for the wide class of non-linear oscillators described ty the Liénard equation.
\begin{enumerate}
\item\label{result1}
A key PMP result is the  bang-bang nature of the optimal driving force: it only takes extreme values, either \(F(t) = -K\) or \(F(t) = +K\). In particular, we have
\begin{equation}
  F(t) = \sgn(p_2(t)) K,  
  \label{eq:F_optimal}
\end{equation}
where \(p_2\) is the second component of the costate \(\bm{p}\). The force changes discontinuously from $F(t)=\pm K$ to $\mp K$ at the times $t_s$ such that $p_2(t_s)=0$, which determine these switching points. Our numerical analysis, detailed below, demonstrates that decreasing the bound 
$K$ increases the number of bangs, leading to a very complex optimal control structure.
\item\label{result2} The final costate \(\bm{p}_f\) is obtained by combining Eqs.~\eqref{eq:hamiltonian_optimal} and \eqref{eq:transversality} at \(t=t_f\). Two cases can be distinguished:
\begin{subequations}\label{eq:final-momenta}
\begin{align}
     \text{(i) }& (x_{1f},x_{2f})\ne (\pm x_{\LC}^{\max},0):  \nonumber \\
    &p_{1f}=\frac{\mu h(x_{1f})x_{2f}+V'(x_{1f})}{F_f \, x_{2f}} ,& p_{2f}&=\frac{1}{F_f}, & p_0&=-1,\label{eq:final-momenta-i}   \\
     \text{(ii) } &(x_{1f},x_{2f})=(\pm x_{\LC}^{\max},0): \nonumber \\ 
    &p_{1f}= \sgn(F_f), & p_{2f}&=0, & p_0&=0, 
    \label{eq:final-momenta-ii}
\end{align}
\end{subequations}
where \((\pm x_{\LC}^{\max},0)\) are the extreme points of the limit cycle over the \(x_1\)-axis.
\item\label{result3} Let us examine optimal trajectories that end at the limit cycle's extreme points---i.e.~case (ii) of result~\ref{result2}, Eq.~\eqref{eq:final-momenta-ii}, which we call critical trajectories.  Over the critical trajectories, all switching points lie on the \(x_1\)-axis: since \(p_0 = 0\), $\forall t$, and at the switching points \(p_2(t_s) = 0\), Eqs.~\eqref{eq:hamiltonian} and \eqref{eq:hamiltonian_optimal} entail \(x_2(t_s)=0\)---PMP excludes the possibility \(p_1(t_s)\neq 0\). Note that the reverse statement is also true: over the critical trajectories, points lying on the $x_1$-axis are switching points.
\item\label{result4} The sign of the force in the final bang can be inferred by analysing the behaviour of \(\chi_{\LC}\) around the final time, as a function of  $\sgn(x_{2f})$ and the position of the initial point with respect to the limit cycle. 
Let us define $x_{2f}^-\equiv x_2(t_f^-)$. For initial points outside the limit cycle, we have: (a) for $x_{2f}^-<0$, the final value of the force is $F_f>0$ and thus the last bang corresponds to $F(t)=+K$, (b) for $x_{2f}^->0$, $F_f<0$ and the last bang corresponds to $F(t)=-K$.
For points inside the limit cycle, the above cases are reversed. 
\end{enumerate}
The analytical results~\ref{result1}--\ref{result4} above allow us to  determine the final conditions of any possible optimal trajectory. These are crucial to obtain a complete numerical solution for the optimisation problem, which we describe below.

\section{Numerical analysis}

Building on the just presented analytical results, we carry out a complete numerical solution of the optimisation problem. Specifically, we present results for the van der Pol oscillator, a paradigmatic case of the Liénard equation with the choice $h(x) = x^2-1$, $V(x) = x^2/2$. The van der Pol equation and some variants have been employed in many different contexts, such as electronic simulation of nervous impulses~\cite{fitzhugh_impulses_1961,nagumo_active_1962}, elastic excitable media~\cite{cartwright_dynamics_1999}, solar cycle~\cite{mininni_stochastic_2000}, active matter~\cite{romanczuk_active_2012}, light-matter interaction~\cite{zhang_van_2024}, or synchronisation of quantum systems~\cite{lee_quantum_2013,ben_arosh_quantum_2021,thomas_quantum_2022,wachtler_topological_2023,impens_shortcut_2023}.

Our numerical approach consists in a backward integration from any final point \(\bm{x}_f^T=(x_{1f},x_{2f})\) over the limit cycle. Equation~\eqref{eq:final-momenta} determines the final value of the momenta and result~\ref{result4} above determines the final bang, which allows us to numerically solve Eq.~\eqref{eq:canon_eq} backward in time from \(t=t_f\). Following Eq.~\eqref{eq:F_optimal}, every time \(p_2\) changes its sign, \(F\) also changes. Starting from all points over the limit cycle, we obtain all the possible trajectories fulfilling  PMP. Finally, to determine the optimal solution for any initial phase space point \(\bm{x}_0^T=(x_{10},x_{20})\), we have two cases:  (i) if there is only one trajectory containing $\bm{x}_0$, then it is the optimal solution; or (ii) if there are multiple trajectories containing $\bm{x}_0$, we choose that reaching the limit cycle the fastest.

To carry out a systematic thorough analysis of the candidates to optimal protocols, we divide the phase plane in two regions: outside the limit cycle, \(\chi_{\LC} > 0 \), and inside the limit cycle, \(\chi_{\LC} < 0\). All numerical results correspond to a damping coefficient $\mu=0.1$, for which the dimensionless relaxation time to the limit cycle is $t_R\simeq 40$~\cite{rios-monje_optimal_2024}.

\subsection{Exterior points: optimal trajectories}

\begin{figure}
  \centering
  \includegraphics[height = 1.77in]{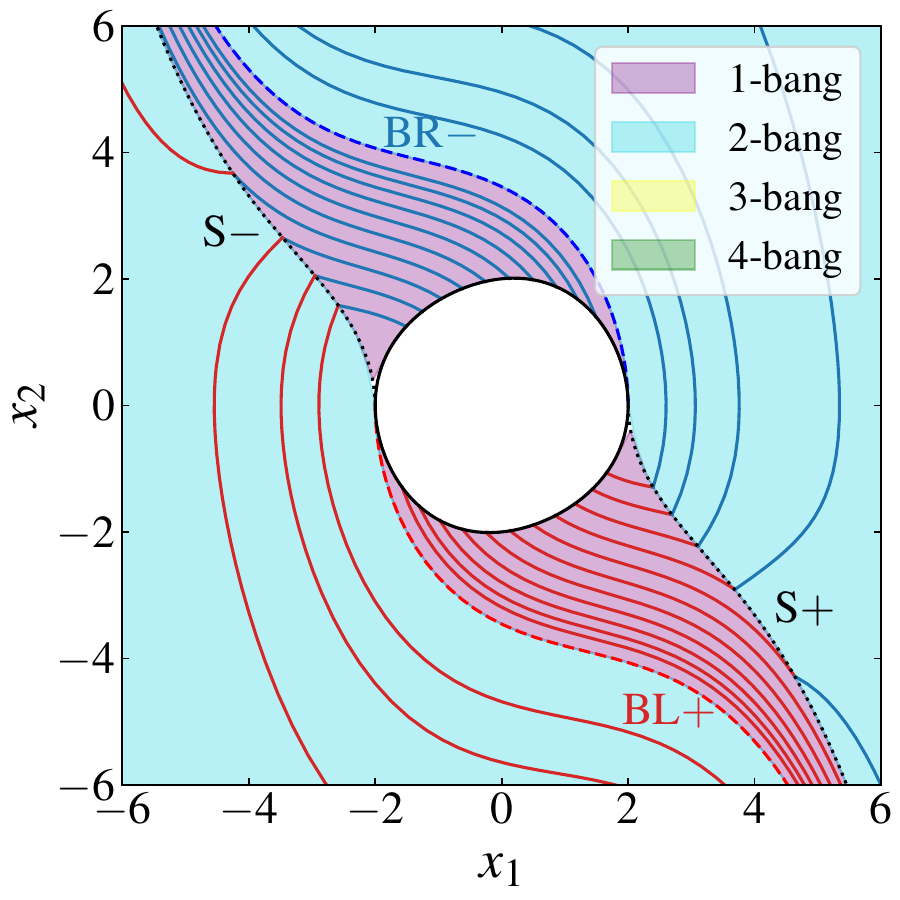}\hspace{-0.5ex}
  \includegraphics[height = 1.77in]{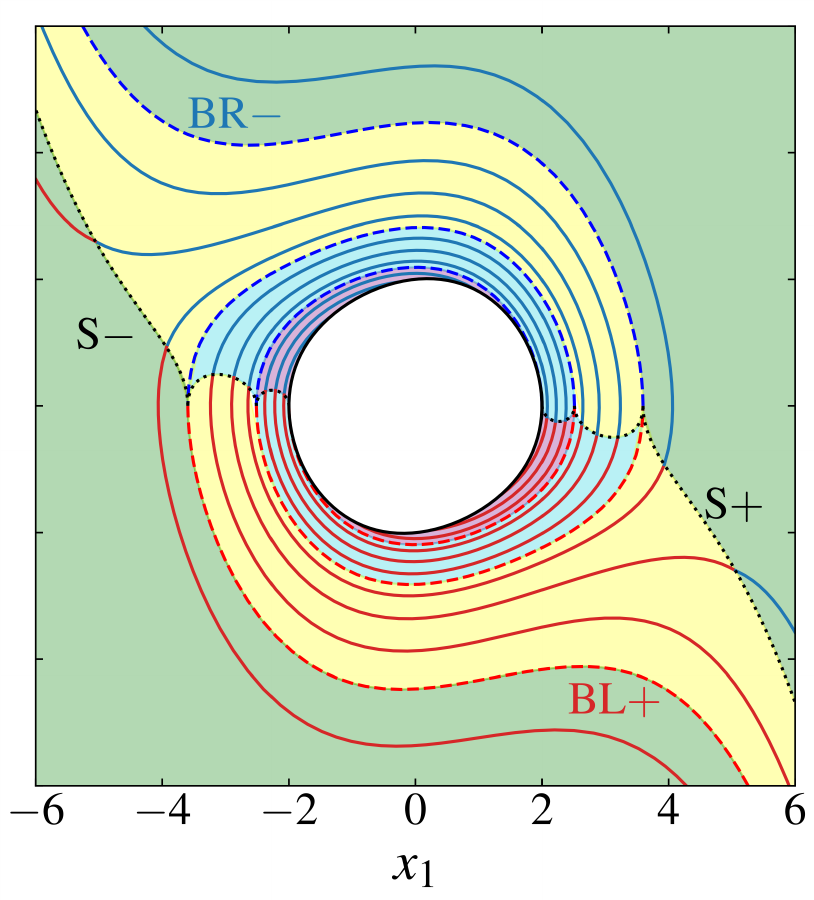}
  \caption{Optimal trajectories for initial phase space points outside the limit cycle. Force bounds are: \(K=2\) (left) and \(K=0.2\) (right), for which the maximum number of bangs is two and four, respectively. The curves delimiting the different regions in phase space are plotted: the switching curves S$+$ and S$-$, at which the force changes sign for $x_1>0$ and $x_1<0$, respectively (dotted); the critical trajectories BL$+$ and BR$-$, obtained by solving the van der Pol equations backward in time from the leftmost $(-x_{\LC}^{\max},0)$ point with $F(0)=+K$ and rightmost point $(+x_{\LC}^{\max},0)$ with $F(0)=-K$, respectively (dashed). Red and blue trajectory arcs correspond to forces $+K$ and $-K$, respectively.}
  \label{fig:outside_k2andk0.2}
\end{figure}

First, we focus on the region outside the limit cycle. Result~\ref{result4} for the final bang, particularised to the van der Pol equation, tells us  $\sgn (F_f) = -\sgn(x_{2f})$  if $x_{2f} \neq 0$, whereas \(\sgn (F_f) = -\sgn (x_{1f}) \) if \(x_{2f} = 0\). With the final sign of the force completely determined, along with the final conditions of the momenta given by Eq.~\eqref{eq:final-momenta}, we apply the above described numerical algorithm for a large enough value of the force bound---results are shown in Fig.~\ref{fig:outside_k2andk0.2}. Since there are no intersections between the plotted trajectories and they are the only ones compatible with PMP, they are the optimal trajectories for any initial point in phase space.

The left panel of Fig.~\ref{fig:outside_k2andk0.2} shows four different regions in phase space, with one-bang or two-bang optimal controls. The regions are delimited by the following curves: the limit cycle, the critical trajectory BL$+$, the critical trajectory BR$-$, and the switching curves S$\pm$, as described in the caption. It is worth stressing that the switching curves S$\pm$ are not trajectories, but the locus of points at which the force changes sign, i.e.~\(p_2 = 0\), over the optimal trajectories. Regarding the optimal control's structure, we have: (i) two bangs with time order $(-K,+K)$, above and to the right of the curves $S+$ and BR$-$; (ii) one bang with $F(t)=+K$ between S$+$ and BL$+$, (iii) one bang with $F(t) = -K$ between S$-$ and BR$-$; (iv) two bangs with  time order $(+K,-K)$ below and to the left of BL$+$ and S$-$. Note that  two bangs are always necessary for initial points on the $x_1$-axis.

The structure above is robust for large enough $K$, such that the critical trajectories BL$+$ and BR$-$ do not intersect the \(x_1\)-axis and thus are one-bang trajectories with no switching points. In the following, we focus the discussion on the behaviour of the BL$+$ trajectory, due to the symmetry of the problem. As $K$ is decreased, there appears a critical value $K_{c1}$ below which BL$+$ intersect the $x_1$ axis once: this intersection  point $x_{c1}(K)$ is a switching point, due to result~\ref{result3}. Prolonging  BL$+$ beyond $x_{c1}(K)$ needs two bangs. If we keep decreasing \(K\), another critical value \(\Kctwo\) appears, below which the trajectory BL$+$ intersect the \(x_1\)-axis twice, at $x_{c1}(K)$ and $-x_{c2}(K)$. Prolonging BL$+$ beyond $-x_{c2}(K)$ needs three bangs, and so on. Note that we have defined all $x_{cn}(K)$ as positive.

\begin{figure}
  \centering
  \includegraphics[width = 3.25in]{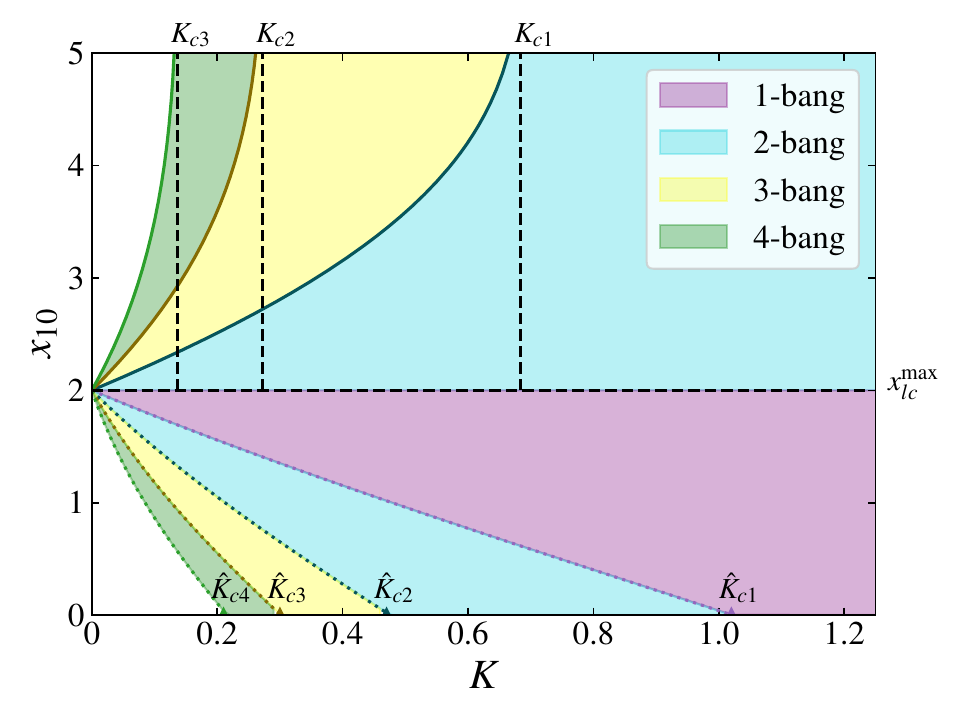}
  \caption{Phase diagram for the number of bangs in the optimal control in the $(K,x_{10})$ plane, for initial points on the \(x_1\)-axis. Color code is the same as in Fig.~\ref{fig:outside_k2andk0.2}. Different regions are delimited by the critical solid (dotted) lines $x_{cn}(K)$ ($\hat{x}_{cn}(K)$) for initial points outside (inside) the limit cycle: to their right, we have $n+1$ ($n$) bangs; to their left, we have $n+2$ ($n+1$) bangs; over them, we have a critical trajectory, either BL$+$ or BR$-$, with \(n\) bangs. Critical lines with $n>4$ are not shown.}
  \label{fig:x0_vs_K_all}
\end{figure}

The above is the signature of the emergence of more complex optimal protocols, with a larger number of bangs, as $K$ is decreased. Indeed, the  curves \(x_{cn}(K)\) are fundamental to understand the complex structure of the optimal control as $K$ is lowered. Figure~\ref{fig:x0_vs_K_all} shows a phase diagram for the number of bangs in the optimal control for initial points on the $x_1$-axis, in which the curves $x_{cn}(K)$ delimit the different regions for $x_{10}>x_{\LC}^{\max}$. 

Summarising, for $K_{c,n+1}<K<K_{cn}$, the critical trajectory BL$+$ intersects the $x_1$ axis at $n$ points and the maximum number of bangs for initial points lying on the $x_1$-axis is $n+2$~\footnote{We are defining $K_{c0}=\infty$ to employ this equation for $n=0$, which is consistent with defining $x_{c0}(K)=x_{\LC}^{\max}$.}. Interestingly, $n+2$ is the maximum number of bangs required to optimally connect any initial phase space point to the limit cycle. The right panel of Fig.~\ref{fig:outside_k2andk0.2} illustrates the specific case $K=0.2$, such that $K_{c3}<K<K_{c2}$, which entails that the maximum number of bangs is four. It is clearly observed that the critical trajectories BL$+$ and BR$-$ intersect the $x_1$-axis at two points: $x_{c1}$ and $-x_{c2}$ for the former, and the symmetric $-x_{c1}$ and $x_{c2}$ for the latter. These points separate different intervals in the $x_1$-axis. For the $x_1$-axis points closest to the limit cycle, $x_{\LC}^{\max} < x_1 < x_{c1}$, two bangs are still enough, as for larger values of $K$. For $x_1=x_{c1}$, the optimal control corresponds to the BL$+$  trajectory and thus have only one bang. For $x_{c1}<x<x_{c2}$, three bangs are necessary.  For $x_1=x_{c2}$, the optimal control corresponds to the BR$-$  trajectory and thus have two bangs: one bang corresponds to the arc between the final point and $(-x_{1c},0)$, and the other to the arc between $(-x_{1c},0)$ and $(x_{2c},0)$. Finally, for the points furthest from the limit cycle, $x_1>x_{2c}$, four bangs are needed.

\subsection{Interior points: optimal trajectories}

For points inside the limit cycle, result~\ref{result4} gives: (i) \(\sgn (F_f)= \sgn( x_{2f}) \) if \(x_{2f} \neq 0\),  \(\sgn (F_f) = \sgn (x_{1f}) \) if \(x_{2f} = 0\). While seemingly simpler due to being a finite region, this case is slightly more complex than that for exterior points: trajectories rewound from different points on the limit cycle may intersect. As a consequence, there may be more than one trajectory starting from an interior point that fulfils PMP. Then, to identify the optimal trajectory, we compare the connection times of all candidate trajectories and select the fastest. Therefore, a coexistence curve emerges, given by the locus of points at which the fastest synchronisation is degenerate.

\begin{figure}
  \centering
   \includegraphics[height = 1.73in]{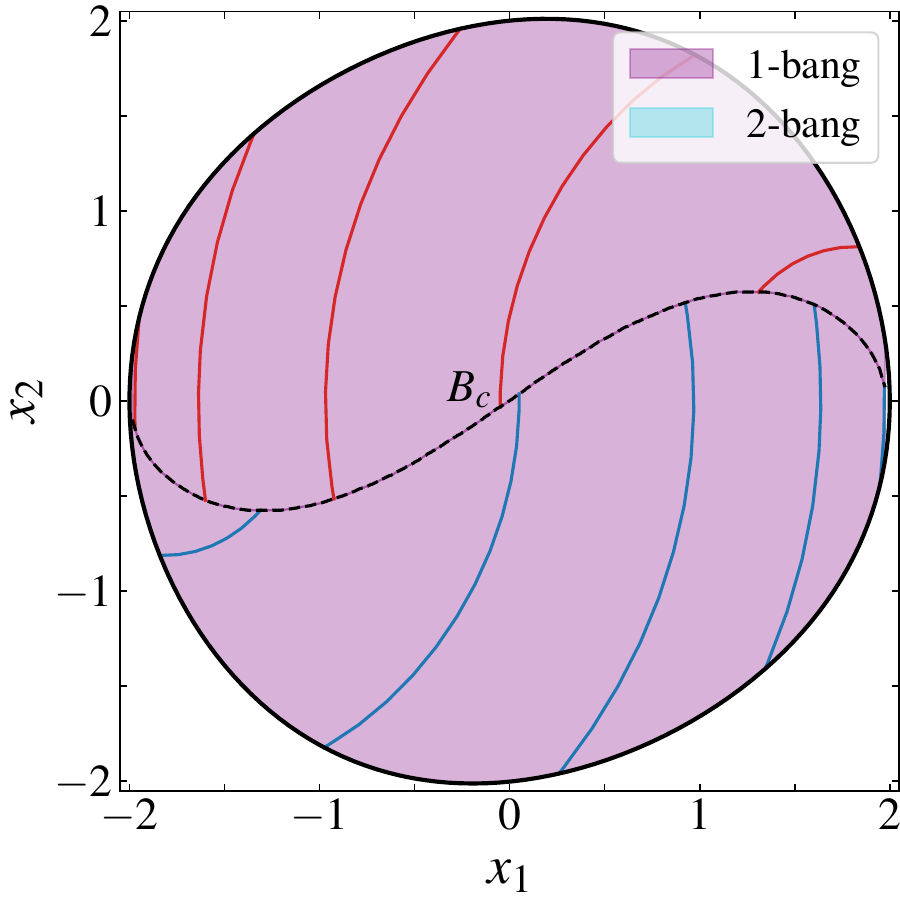}
  \includegraphics[height = 1.73in]{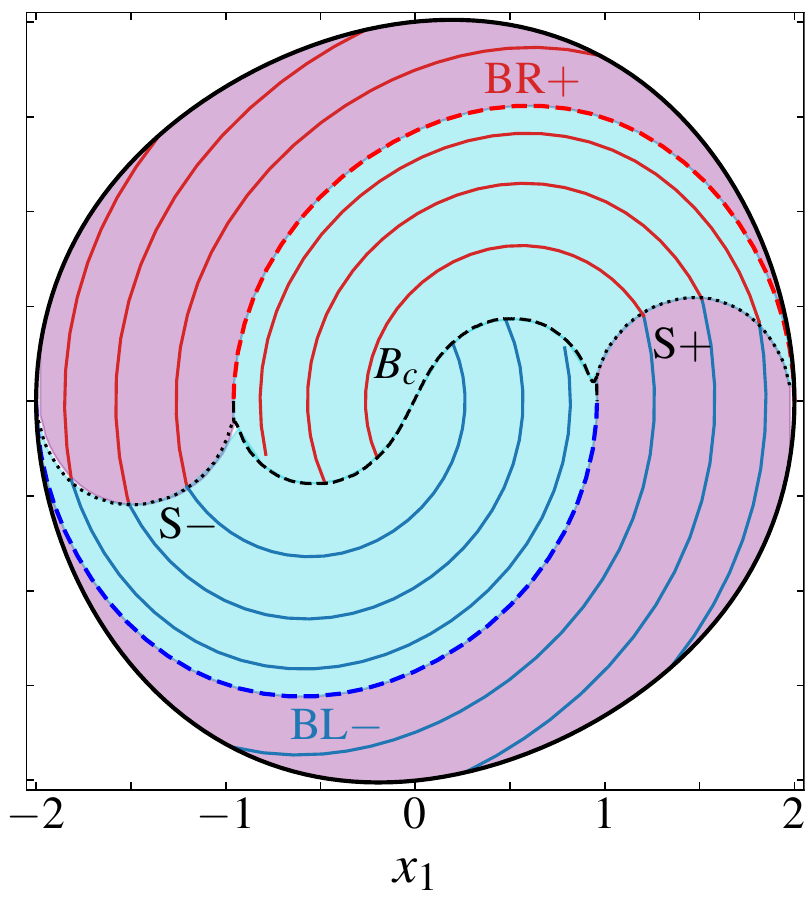}
  \caption{Optimal trajectories for initial phase space points inside the limit cycle. Force bounds are $K=2$ (left) and $K=0.5$ (right). The number of bangs required to connect any initial point to the limit cycle is: (left) always one; (right) one or two, depending on the initial point. In addition to the switching curves S$\pm$ and the critical trajectories BL$+$ and BR$-$, analogous to those in Fig.~\ref{fig:outside_k2andk0.2}, we have plotted the coexistence curve $B_c$ (black dashed), at which the intersecting one-bang trajectories give the same connection time. 
  }
  \label{fig:inside_k2_and_05}
\end{figure}

For $x_{10}<x_{\LC}^{\max}$, Fig.~\ref{fig:x0_vs_K_all} illustrates the critical curves \(\hat{x}_{cn}(K)\) governing the optimal control structure inside the limit cycle. Analogously to the case of exterior points, these critical curves are determined by the intersection points with the \(x_1\)-axis of the inside critical trajectory ending at \((-x_{\LC}^{\max},0)\),  as a function of $K$. In Fig.~\ref{fig:inside_k2_and_05}, we observe that (i) one bang is enough for any interior point for \(K=2\) (left panel) (ii) either one or two bangs are necessary for $K=0.5$ (right panel), consistently with the phase diagram in Fig.~\ref{fig:x0_vs_K_all}. Smaller values of $K$ would involve an increasingly larger number of bangs. 

\subsection{Optimal syncrhonisation time}

\begin{figure}
  \centering
  \includegraphics[height = 1.4in]{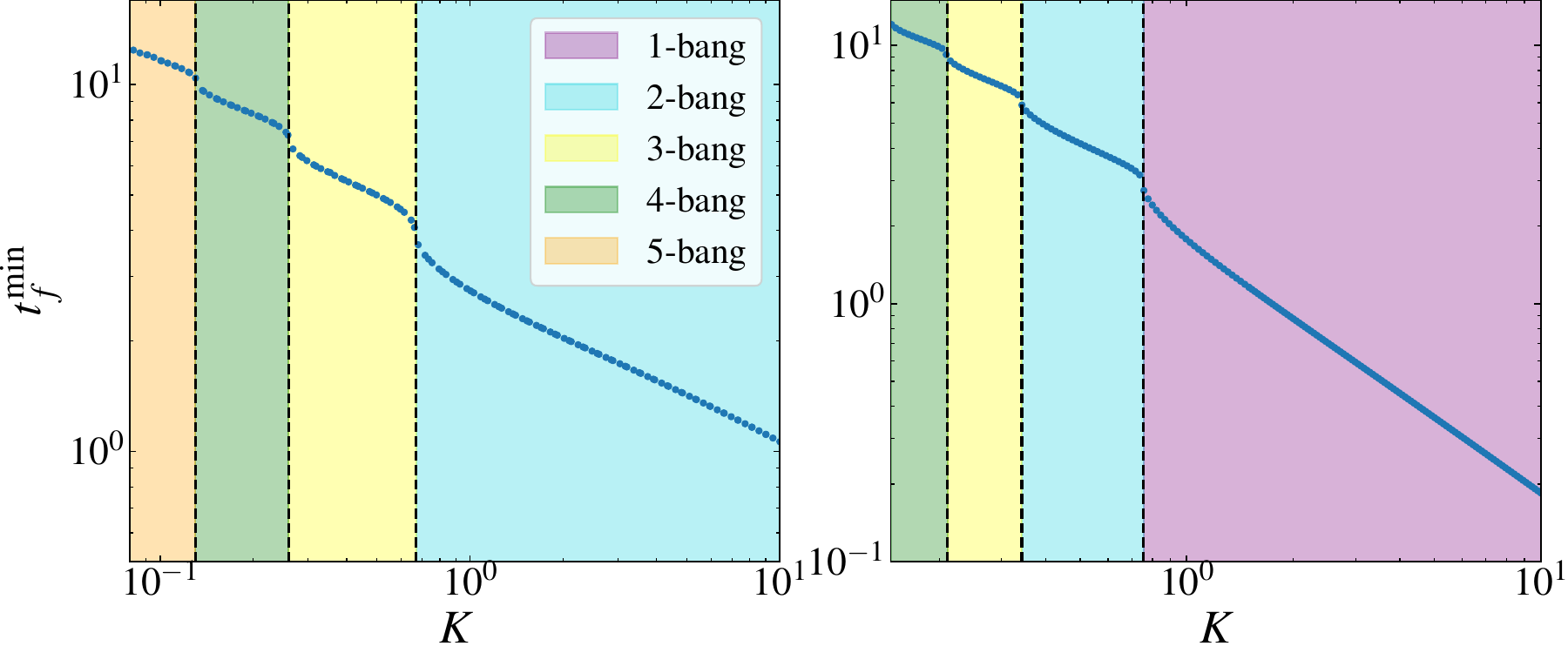}
  \caption{Minimum connection time as a function of the force bound. Specifically, we have chosen two initial points on the $x_1$-axis, one exterior $(5,0)$ (left) and one interior  $(1,0)$ (right). In both panels, the vertical lines correspond to the values $K^*$ such that $x_{cn}(K^*)=x_{10}$ (left) and $\hat x_{cn}(K^*)=x_{10}$ (right). The minimum connection time is continuous at these critical values $K^*$, but its first derivative with respect to $K$ is not. 
  }
  \label{fig:tfmin}
\end{figure}

We present the minimum connection time as a function of the force in Fig.~\ref{fig:tfmin}. As the bound in the force $K$ decreases, the connection time consistently increases. The discontinuous derivative of the minimum connection time at the critical values $K_{cn}$ (left, exterior point) and $\hat{K}_{cn}$ (right, interior point) is a consequence of the change in the number of bangs of the optimal protocol thereat, which was illustrated in Fig.~\ref{fig:x0_vs_K_all}. Inside each $n$-bang region and for points not very close to the critical values $K_{cn}$ and $\hat{K}_{cn}$, the behaviour $t_f^{\min}$ as a function of $K$ is compatible with a power law; $\ln t_f^{\min}$ displays an approximate linear behaviour with $\ln K$.

\section{Conclusions}

In this paper, we have solved the minimum-time synchronisation problem to the Liénard's limit cycle. Leveraging Pontryagin's Maximum Principle, we have derived general analytical results, \ref{result1}--\ref{result4}, that have been instrumental for developing a robust backward-integration algorithm for finding numerical solutions.

Our application of such general framework to the paradigmatic case of the van der Pol oscillator has revealed detailed insights of the optimal control structure across the whole phase space. As the force bound $K$ decreases or the initial condition moves further from the limit cycle, both the complexity of the optimal control and the minimum synchronisation time increase. This behaviour is intuitively sound, as greater constraints on the system lead to more intricate trajectories and longer connection times.

Crucial for our analysis is the key role played by ``critical trajectories,'' i.e.~those rewound from the extreme points of the limit cycle. These trajectories define critical $K$ values, both for exterior ($K_{cn}$) and interior ($\hat{K}_{cn}$) points. For a given force bound $K$, these critical values ultimately determine the maximum number of control bangs required for the optimal time connection to the limit cycle from any initial phase space point.

The results presented here pave the way for future research in synchronisation and self-sustained oscillations. Further extension of our analytical results, e.g.~by finding analytical expressions for the critical values of $K$ and the corresponding critical curves \(x_{cn}(K)\) and \(\hat{x}_{cn}(K)\), is an interesting prospect. Also, given the relevance of synchronisation across numerous fields, applying the ideas introduced here to more complex and realistic models is another promising avenue for exploration.

\section*{Acknowledgments}

This publication has been funded  by  Junta de Andalucía-Consejería de Universidad, Investigación e Innovación, in the frame of the research project with reference ProyExcel\_00796. CRM, CAP and AP acknowledge support from Grant PID2024-155268NB-I00 funded by MICIU/AEI/10.13039/501100011033/ FEDER, UE, and also from
the applied research and innovation Project PPIT2024-31833, cofunded by EU--Ministerio de Hacienda y Función Pública--Fondos Europeos--Junta de Andalucía--Consejería de Universidad, Investigación e Innovación. DGO acknowledges support from Institut Universitaire de France. CRM acknowledges support from the Spanish Ministry of Science, Innovation and Universities FPU programme through Grant FPU22/01151.

\section*{Data Availability Statement}
The Python codes employed for generating the data are openly available in the~\href{https://github.com/fine-group-us/minimum_time_synchronisation}{GitHub page} of University of Sevilla's FINE research group.

\appendix

\section{Derivation of results \ref{result1}--\ref{result4} from PMP} 

In this section, we detail the proofs of the analytical results \ref{result1}--\ref{result4}.

\subsection{Result 1: the optimal control is bang-bang}

To prove this, we use property (ii) of PMP, which states that the  Hamiltonian is maximum with respect the control at the optimal one. Since the control set is bounded, \(F(t)\in[-K,K]\), the maximum of \(\mathcal{H}\) may lie either (i) in the open interval \((-K,K)\) or (ii) at the boundaries of the control set, i.e.~\(F(t) = \pm K\). In the first case, bringing to bear Eq.~\eqref{eq:hamiltonian}, we get
\begin{equation}
    \partial_F \mathcal{H} = 0 \Rightarrow p_2 = 0.
\label{eq:dHdF}
\end{equation}

Due to the linear dependence of the Hamiltonian on \(F\), this is a singular optimal control problem and Eq.~\eqref{eq:dHdF} does not give us any information about the optimal control. Therefore, we apply the Kelley condition~\cite{leitmann_topics_1967,zelikin_singular_2005} to get further information about the control. We start by differentiating Eq.~\eqref{eq:dHdF} with respect to time:
\begin{equation}
  \frac{d}{dt}\left(\partial_F\mathcal{H}\right) = 0 \Rightarrow \dot{p}_2 = 0.
\label{eq:dtdHdF}
\end{equation}
Taking into account \(p_2\) equation in Eq.~\eqref{eq:canon_eq}, we get 
\begin{equation}
  0 = - p_1 + \mu h(x_1)\cancelto{0}{p_2} = -p_1.
\label{eq:p1_constant}
\end{equation}
Therefore, if there is a time interval \((t_1,t_2)\) where Eq.~\eqref{eq:dHdF} is satisfied, we conclude that \(\bm{p} = \bm{0}\) therein. That would mean, using Eq.~\eqref{eq:hamiltonian_optimal}, that also \(p_0 = 0\) and \(\{\bm{p},p_0\} = \{\bm{0},0\}\) for all \(t\in(t_1,t_2)\). But this is incompatible with PMP and case (i) has to be discarded. Therefore, the only scenario is (ii); over the possible optimal control we have \(F(t) = \pm K\). The sign of the control is determined by the sign of the coefficient of \(F\) in \(\mathcal{H}\), which is \(p_2\), which leads to Eq.~\eqref{eq:F_optimal}.

\subsection{Result 2: final value of the costate}

Let us know focus on the transversality condition, Eq.~\eqref{eq:transversality}, and Eq.~\eqref{eq:hamiltonian_optimal} at the final time \(t_f\),
\begin{subequations}
  \begin{align}
    p_{1f}x_{2f}+p_{2f}\left[-\mu h(x_{1f})x_{2f}-V'(x_{1f})\right]&=0, \\
    p_{1f}x_{2f} + p_{2f}\left[-\mu h(x_{1f})x_{2f}-V'(x_{1f})+F_f\right] + p_0 &= 0,
  \end{align}
\label{eq:final_cond}
\end{subequations}
respectively, where \(\bm{p}^T(t_f)=(p_{1f},p_{2f}) \) and \(F_f \equiv F(t_f)\). Combining these two equations yields
\begin{equation}
  p_{2f}F_f + p_0 = 0.
\label{eq:final_cond_p2F} 
\end{equation}
The limit cycle of the Liénard equation has two extreme points $(\pm x_{\LC}^{\max},0)$, at which $x_1$ is either maximum or minimum and thus $\dot{x}_1=x_2$ vanishes. It is only at these points of the limit cycle that the vector $\bm{d}$ of its tangent space is vertical and $x_2$ vanishes. On the one hand, if the final point $(x_{1f},x_{2f})$ is not one of the extreme points, then we have that $p_{2f}\ne 0$---because the final vector of momenta is not horizontal---and Eq.~\eqref{eq:final_cond_p2F} entails that $p_0\ne 0$. Normalising it to $p_0=-1$, as usual in optimal control theory~\cite{pontryagin_mathematical_1987,liberzon_calculus_2012}, we obtain Eq.~\eqref{eq:final-momenta-i} in result~\ref{result2}.  On the other hand, if the final point $(x_{1f},x_{2f})$ is one of the extreme points $(\pm x_{\LC}^{\max},0)$, then we have that $p_{2f}= 0$ and Eq.~\eqref{eq:final_cond_p2F} now implies $p_0=0$: we need $p_{1f}\ne 0$ to fulfil PMP. The sign of $p_{1f}$ is determined by the sign of $F_f$. Considering Eq.~\eqref{eq:canon_eq} for $p_2$, with $p_{2f}=0$, we get $\dot{p}_{2f} = -p_{1f}$.
Now if we assume $F_f > 0$, then, by Eq.~\eqref{eq:F_optimal}, $p_{2f}^{-}>0$. Since $p_{2f}=0$, this implies $\dot{p}_{2f}< 0 $ that in turn entails $p_{1f}>0$. Conversely, if $F_f<0$, we get the opposite result. This gives Eq.~\eqref{eq:final-momenta-ii} in result~\ref{result2}.

\subsection{Result 4: sign of the final  bang}

Bringing to bear the final value of the force \(F(t_f)=F_f\), we can rewind the evolution equations from the final point over the limit cycle \(\bm{x}_f^T=(x_{1f},x_{2f})\). We have taken \(\chi_{\LC}>0 \) ($<0$) for phase plane points lying outside (inside) the limit cycle. Therefore,  $\nabla\chi_{\LC}(x_{1f},x_{2f})$ points outwards from limit cycle, and since the gradient is perpendicular to the tangent to the limit cycle, we write
\begin{align}
    \nabla\chi_{\LC}(\bm{x}_f)^T&=c_f(\mu h(x_{1f})x_{2f}+V'(x_{1f}),x_{2f}), \; c_f>0.
    \label{eq:grad_chi}
\end{align}
Note that $c_f$ may depend on  $\bm{x}_f$, $c_f=c(\bm{x}_f)$.

Now we consider an infinitesimal time interval \((t_f-\varepsilon,t_f)\) and the corresponding phase space point \(\bm{x}_f^-=(x_{1f}^-,x_{2f}^-)\) at \(t_f-\varepsilon\). Up to first order in \(\varepsilon\), we have
\begin{subequations}
\begin{align}
x_{1f}^-&=x_{1f}-\varepsilon \,x_{2f},\\
x_{2f}^-&=x_{2f}-\varepsilon \left[-\mu h(x_{1f})x_{2f}-V'(x_{1f})+F_f\right], 
\label{eq:x_minus_b}
\end{align}  
\label{eq:x_minus}
\end{subequations}
in which we have neglected $O(\varepsilon^2)$ terms. 
Using Eqs.~\eqref{eq:x_minus} to expand \(\chi_{\LC}\) at \(\bm{x}_f\) yields
\begin{equation}
    \chi_{\LC}(\bm{x}_f^-)=-\varepsilon F_f \left.\partial_{x_2}\chi_{\LC}\right|_{\bm{x}_f}=-\varepsilon c_f x_{2f} F_{f},
    \label{eq:chi_firstorder}
\end{equation}
where we have omitted O(\(\varepsilon^2\)) terms. For $x_{2f}\ne 0$, Eq.~\eqref{eq:chi_firstorder} entails that $x_{2f} F_f<0$ when reaching the limit cycle from the outside (exterior initial point), whereas $x_{2f} F_f>0$ when reaching the limit cycle from the inside (interior initial point)---in agreement with result~\ref{result4}.

We have to study separately the case \(x_{2f}=0\), in which the first order approximation in  Eq.~\eqref{eq:chi_firstorder} vanishes. This corresponds to the critical trajectories in which $x_{1f}=\pm x_{\LC}^{\max}$. Expanding \(\chi_{\LC}\) from these extreme points of the limit cycle to second order in \(\varepsilon\) gives
\begin{equation}
    \chi_{\LC}(\bm{x}_f^-)=-\varepsilon c_f \,x_{2f}^-\, F_f/2,
\end{equation}
once more using Eq.~\eqref{eq:x_minus}. For initial points outside the limit cycle, for which $\chi_{\LC}>0$, we get: (i) if $x_{2f}^-<0$, then $F_f >0$ and the last bang must be $F(t) = +K$; (ii) if $x_{2f}=0$ but $x_{2f}^->0$, then $F_f <0$ and the last bang must be $F(t) = -K$. For initial points inside the limit cycle, the sign of the bangs is reversed, taking into account that $\chi_{\LC}<0$. This completes the proof of our result~\ref{result4}.

Equation~\eqref{eq:x_minus_b} leads to
\begin{equation}
    x_{2f}=0 \implies x_{2f}^-=\varepsilon \left[V'(\pm x_{\LC}^{\max})-F_f\right].
\end{equation}
Therefore, for small values of the bound, such that $V'(x_{\LC}^{\max})>K$, the sign of $x_{2f}^-$ is perfectly determined: negative at $(-x_{\LC}^{\max},0)$ and positive at $(x_{\LC}^{\max},0)$. This leads to the critical trajectories BL$+$ and BR$-$ in the main text. For larger values of the bound, $V'(x_{\LC}^{\max})<K$, $x_{2f}^-$ may have both signs and there may exist critical trajectories BL$-$ (with $x_{2f}^->0$) and BR$+$ (with $x_{2f}^-<0$), in addition to the trajectories BL$+$ and BR$-$. However, the trajectories  BL$-$ and BR$+$ consist of points that belong to other candidates to optimal trajectories with shorter synchronisation times, and therefore they can be discarded.

\bibliography{bibliography}

\end{document}